\documentclass{SCAE}
\usepackage{amsmath,amsfonts,amsthm,amssymb,amscd,epsfig}
\usepackage{graphicx,psfrag,overpic}
\usepackage{yhmath}
\usepackage[mathscr]{eucal}
\usepackage{xspace}

\newcommand{\R}{\mathbb{R}}

\newcommand{\CS}{\mathcal{S}}

\newcommand{\xx}{\mathbf{x}}

\newcommand{\ff}{\mathbf{f}}
\newcommand{\kk}{\mathbf{k}}

\newcommand{\uu}{\mathbf{u}}

\newcommand{\ga}{\alpha}
\newcommand{\gb}{\mathbf{\beta}}

\newcommand{\gD}{\Delta}
\newcommand{\po}{\partial}

\newcommand{\ve}{\varepsilon}
\newcommand{\vae}{\varepsilon}

\newcommand{\gd}{\delta}

\newcommand{\del}{{\partial}}
\renewcommand{\O}{{O}}
\renewcommand{\d}{\delta}

\renewcommand{\b}{\beta}

\newcommand{\s}{\sigma}
\newcommand{\g}{\gamma}

\newcommand{\Om}{\Omega}

\newcommand{\grad}{\nabla}

\newcommand{\bea}{\begin{eqnarray}}
\newcommand{\eea}{\end{eqnarray}}

\newcommand{\yy}{\mathbf{y}}

\newcommand{\pw}{\partial \mathcal{W}}

\newcommand{\e}{\varepsilon}
\newcommand{\MM}{\mathbb{M}}
\newcommand \der{\partial}
\newcommand \vphi{\varphi}

\newcommand \gam{\gamma}
\newcommand \alp{\alpha}

\newcommand \ol{\overline}

\newcommand{\nnu}{{\boldsymbol\nu}}

\newcommand \Wedge{\Gamma_{\rm wedge}}

\newcommand \shock{\Gamma_{\rm shock}}

\newcommand \lefttop{P_1}
\newcommand \righttop{P_2}

\newcommand \leftbottom{P_4}
\newcommand \rightbottom{P_3}

\newcommand \ivphi{\varphi_0}
\newcommand \iu{u_{10}}

\newcommand \irho{\rho_0}

\newcommand \leftshockop{\mathcal{S}_{0}}
\newcommand \leftsonicop{\Gamma_{\rm sonic}^{1}}

\newcommand \oOmop{\Om_1}

\newcommand \rightshockop{\mathcal{S}_{1}}
\newcommand \rightsonicop{\Gamma_{\rm sonic}^{2}}

\newcommand \nOmop{\Om_2}
\newcommand \bPhi{\bar{\Phi}}
\newcommand{\xxi}{{\boldsymbol\xi}}

\numberwithin{equation}{section}
\numberwithin{figure}{section}
\begin{document}

%Basic Information
\Year{2017} %
\Month{May}
%{January}
\Vol{60} %
\No{2} 
\BeginPage{1} %
\EndPage{19} %
\AuthorMark{Chen G -Q}
\ReceivedDay{December 21, 2016}
%\AcceptedDay{February 22, 2017}
%\PublishedOnlineDay{; published online May 6, 2017}
\DOI{10.1007/s11425-000-0000-0} % The author doesn't need fill in it.

% \title[short text for running head]{full title}{comments for title}
\title[Supersonic Flow onto Solid Wedges, Multidimensional Shock Waves and Free Boundary Problems]
%[M-D Shock Waves and Free Boundary Problems]
{Supersonic Flow onto Solid Wedges,\\ Multidimensional Shock Waves \\ and Free Boundary Problems}
{Dedicated to the 80th Birthday of Professor Tatsien Li}
%\title[SCIENCE CHINA Mathematics  journal sample]{SCIENCE CHINA Mathematics  journal sample}{}

% \author[]{Full name}{footnote}
% Remark:  One \author for one author

\author[1,2,3,4]{CHEN  Gui-Qiang}{}
%\author[3,4]{LASTNAME2 FirstName2}{Corresponding author}
%\author[4]{LASTNAME3 FirstName3}{}

%
\address[{\rm 1}]{School of Mathematical Sciences, Fudan University, Shanghai 200433, P. R. China;}
\address[{\rm 2}]{Academy of Mathematics and Systems Science, Chinese Academy of Sciences, Beijing 100190, P.
R. China;}
\address[{\rm 3}]{University of Chinese Academy of Sciences, Beijing 100049, P. R. China;}

\address[{\rm 4}]{Mathematical Institute, University of Oxford, Oxford, OX2 6GG, UK}
%\address[{\rm2}]{Department of Mathematics, University2, City2 {\rm 100002}, Country2;}
%\address[{\rm3}]{Department of Mathematics, University3, City3 {\rm100003}, Country3;}
%\address[{\rm4}]{College of Science, University4, City4 {\rm100004}, Country4}
\Emails{$\,$  chengq@maths.ox.ac.uk}
% first@mail.com,second@mail.com, third@mail.com}
\maketitle

%     Abstract is required.

{\begin{center}
\parbox{14.5cm}
{\begin{abstract}
When an upstream steady uniform supersonic flow
impinges onto a symmetric straight-sided wedge, governed by the Euler equations,
there are two possible
steady oblique shock configurations if the wedge angle is less than
the detachment angle -- the steady weak shock with supersonic or subsonic
downstream flow
(determined by the wedge angle that is less or larger than
the sonic angle)
and the steady strong shock
with subsonic downstream flow, both of which satisfy the entropy
condition.
The fundamental issue -- whether one or both of the steady weak and strong shocks are physically
admissible solutions -- has been vigorously debated over the past eight decades.
In this paper, we survey some recent developments on the stability analysis of the steady shock solutions
in both the steady and dynamic regimes.
For the static stability, we first show how the stability problem can be formulated as an initial-boundary value type problem
and
then reformulate it into a free boundary problem when the perturbation of
both the upstream steady supersonic flow and the wedge boundary
are suitably regular and small, and we finally present some recent results on the static stability of the steady supersonic and
transonic shocks.
For the dynamic stability for potential flow,
we first show how the stability problem can be formulated as an initial-boundary value problem
and
then use the self-similarity of the problem to reduce it into a boundary value problem and further reformulate it into a
free boundary problem,
and we finally survey some recent developments in solving this free boundary problem for the existence of the Prandtl-Meyer
configurations
that tend to the steady weak supersonic or transonic oblique shock solutions as time goes to infinity.
Some further developments and mathematical challenges in this direction are also discussed.
\vspace{-3mm}
\end{abstract}}

\end{center}}

%  Keyword is required.

%  \subjclass is required.

%%%%%%%%%%%%%%%%%%%%%%%%%%%%%%%%%%%%%%%%%%%%%%%%%%%%%%%%%%%%
\renewcommand{\baselinestretch}{1.2}
\begin{center} \renewcommand{\arraystretch}{1.5}
{\begin{tabular}{lp{0.8\textwidth}} \hline \scriptsize
{\bf Citation:}\!\!\!\!&\scriptsize CHEN Gui-Qiang. \makeatletter\@titlehead.
Sci China Math, 2016, 59,
 %\@Year, \@Vol: \@BeginPage--\@EndPage,
 doi:~\@DOI\makeatother\vspace{1mm}
\\
\hline
\end{tabular}}
\end{center}

%%%%%%%%%%%%%%%%%%%%%%%%%%%%%%%%%%%%%%%%%%%%%%%%%%%%%%%%%%%%
%% Text of article.
%%%%%%%%%%%%%%%%%%%%%%%%%%%%%%%%%%%%%%%%%%%%%%%%%%%%%%%%%%%%
%    Section headings
\baselineskip 11pt\parindent=10.8pt  \wuhao

\section{Introduction}

We survey some recent developments in the analysis
of supersonic flow onto solid wedges (see Fig. 1.1),
involving multidimensional shock waves, and related initial-boundary value type
problems and free boundary problems
for the Euler equations for compressible fluids.
This paper is dedicated to Professor Tatsien on the occasion of his 80th birthday,
who has made pioneering and fundamental contributions to this research direction
and related areas ({\it cf}.
\cite{GuLi, Li-0, Li, LiYu1, LiYu2, LiYu3}, \cite{GuLiHou1,GuLiHou2,GuLiHou3},
and the references cited therein), as we discussed below.

The wedge problem is a longstanding fundamental problem in
mathematical fluid mechanics, partly owing to both the rich
wave configurations
in the fluid flow around the wedge
and the mathematical challenges involved.
More importantly, the solution configurations of the wedge problem
are core configurations in the structure of
global steady entropy solutions,
as well as global dynamic entropy solutions
of the two-dimensional Riemann problem,
for multidimensional hyperbolic systems
of conservation laws.
The steady and Riemann solutions themselves are expected to be local building blocks
and determine local structures, global attractors, and large-time asymptotic
states of general entropy solutions for the systems.
In this sense, we have to understand the solution
configurations and their stability in order to understand fully the global
entropy solutions of the multidimensional
hyperbolic systems of conservation laws.

\begin{figure}
 	\centering
 	\includegraphics[height=35mm]{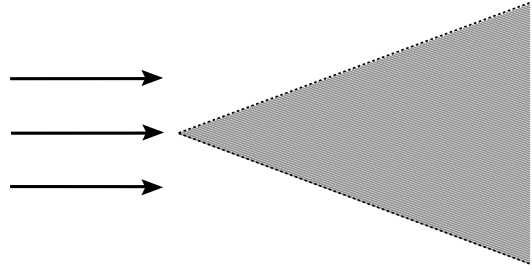}
		\caption{Upstream steady uniform supersonic flow past a symmetric straight-sided wedge}
\label{fig-1}
 \end{figure}

The two-dimensional steady, full Euler equations take the form:
 \begin{equation}\label{Euler1}
 \left\{\begin{aligned}
 & \nabla_\xx\cdot (\rho \uu)=0, \\
 &\nabla_\xx\cdot\left(\rho{\uu\otimes\uu}\right)
 +\nabla_\xx p=0,\\
 &\nabla_\xx\cdot\big(\rho\uu(E+\frac{p}{\rho})\big)=0,
 \end{aligned}\right.
  \end{equation}
 where $\grad_\xx$ is the gradient in $\xx=(x_1,x_2)\in\R^2$,
 $\uu=(u_1,  u_2)$ the velocity, $\rho$ the density, $p$ the pressure,
 as well as
 \begin{equation}\label{1.1a}
 E=\frac{1}{2}|\uu|^2+ e
 \end{equation}
 is the total energy with the internal energy $e$.
 The other two thermodynamic variables are temperature $T$ and
 entropy $S$. If $\rho$ and $S$ are chosen as the independent variables,
 then the constitutive relations can be written as
 $$
 (e, p, T)=(e(\rho, S), p(\rho, S), T(\rho, S))
 $$
 governed by
 $
 TdS=de -\frac{p}{\rho^2}d\rho.
 $
 For an ideal gas,
\begin{equation}\label{1.2a}
 p=R\rho T, \qquad  e=c_vT, \qquad\gamma=1+\frac{R}{c_v}>1,
 \end{equation}
 and
 \begin{equation}\label{1.3a}
  p=p(\rho, S)=\kappa \rho^\gamma {\rm e}^{S/c_v}, \quad
 e=\frac{\kappa}{\gamma-1}\rho^{\gamma-1}{\rm e}^{S/c_v}=\frac{RT}{\gamma-1},
 \end{equation}
 where $R,\kappa$, and $c_v$ are all positive constants.

The sonic speed
of the polytropic gas flow is
\begin{equation}\label{1.4a}
c=\sqrt{\frac{\gamma p}{\rho}}.
\end{equation}
The flow is subsonic if $|\uu| < c$ and supersonic
if $|\uu|> c$.
For a transonic flow, both cases occur in the flow, and then system \eqref{Euler1}
is of mixed-composite elliptic-hyperbolic type, which consists
of two equations of mixed elliptic-hyperbolic type
and two equations of transport-type (which are hyperbolic).

System \eqref{Euler1} is a prototype of general nonlinear systems
of conservation laws:
\begin{equation}\label{Euler1a}
\nabla_\xx\cdot \mathbf{F}(U)=0,  \qquad\,\, \xx\in \R^n,
\end{equation}
where $U: \R^n\to \R^m$ is unknown, while $\mathbf{F}: \R^m\to \MM^{m\times n}$
is a given nonlinear mapping for the $m\times n$ matrix space $\MM^{m\times n}$.
For \eqref{Euler1}, we may choose $U=(\uu, p, \rho)$.
The systems with form \eqref{Euler1a}
often govern time-independent solutions of multidimensional quasilinear
hyperbolic systems of conservation laws; {\it cf.} Dafermos \cite{Dafermos}
and Lax \cite{Lax}.

It is well known that, for an upstream steady uniform supersonic
flow past a symmetric straight-sided wedge (see Fig. 1.1):
\begin{equation}\label{wedge-1}
W:=\{\xx=(x_1,x_2)\in \R^2\,:\, |x_2|<x_1\tan \theta_{\rm w}, x_1>0\}
\end{equation}
whose angle $\theta_{\rm w}$
is less than
the detachment angle $\theta^{\rm d}_{\rm w}$,
there exists an oblique shock emanating from the
wedge vertex.
Since the upper and lower subsonic regions do not interact with each other,
it suffices to study the upper part.
Then, more precisely, if the upstream steady flow is a uniform supersonic state,
we can find the corresponding
constant downstream flow along the straight-sided upper wedge boundary,
together with a straight shock separating
the two states.
The downstream flow is determined by the shock polar
whose states in the phase space are governed by
the Rankine-Hugoniot conditions
and the entropy condition; see Fig.~1.2 and \S 2.
Indeed, Prandtl in \cite{Prandtl} first employed the shock polar analysis
to show that there are two possible
steady oblique shock configurations when the wedge angle $\theta_{\rm w}$ is less than
the detachment angle $\theta_{\rm w}^{\rm d}$ -- The steady weak shock with supersonic or subsonic downstream flow
(determined by the wedge angle that is less or larger than
the sonic angle $\theta^{\rm s}_{\rm w}$)
and the steady strong shock
with subsonic downstream flow, both of which satisfy the entropy
condition, provided that no additional conditions are assigned at downstream.
See also Busemann \cite{Busemann}, Courant-Friedrichs \cite{CourantF}, Meyer \cite{Meyer},
and the references cited therein.

\begin{figure}
 	\centering
 	\includegraphics[height=35mm]{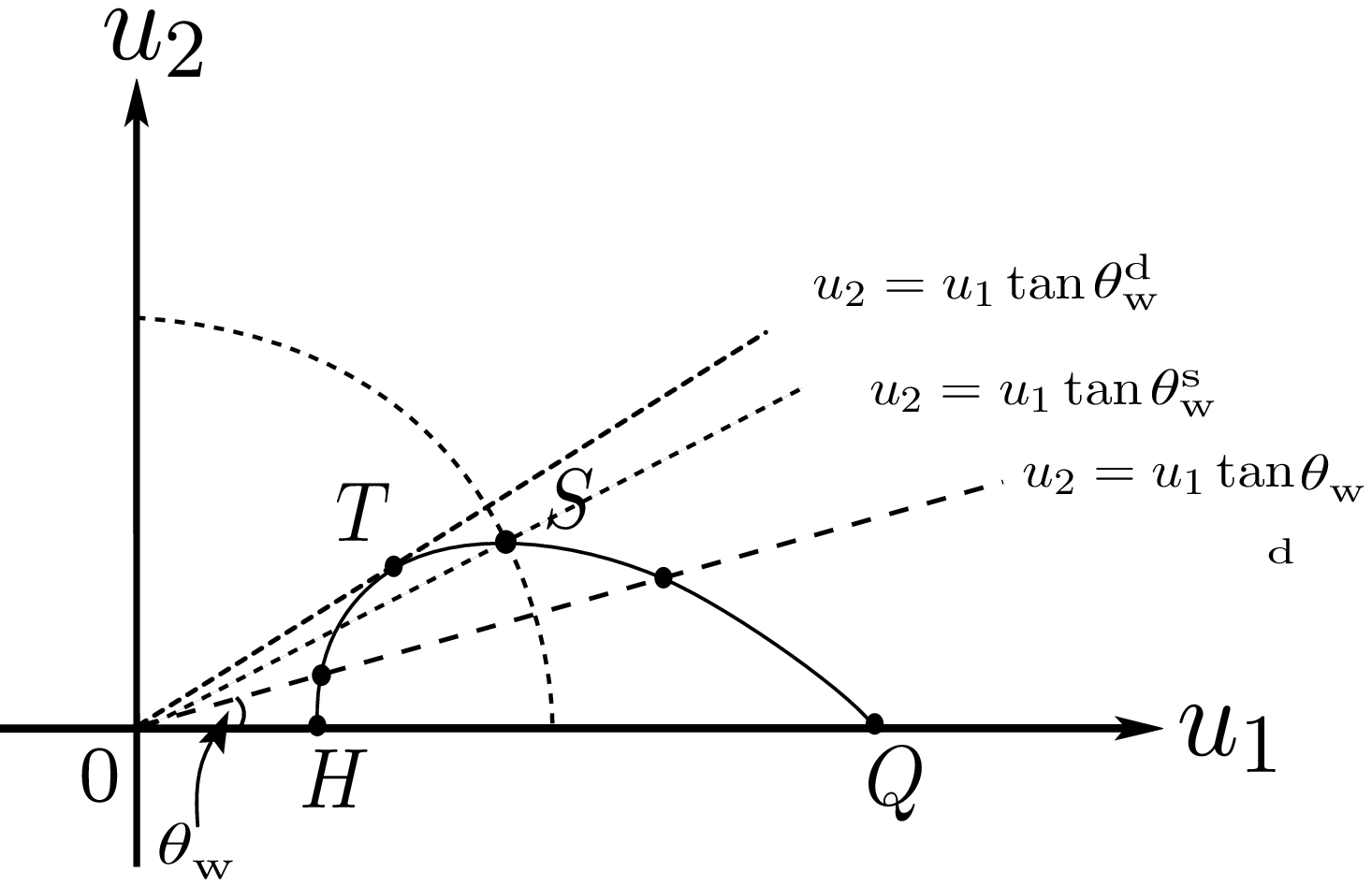}
\hspace{12mm}
 		\includegraphics[height=35mm]{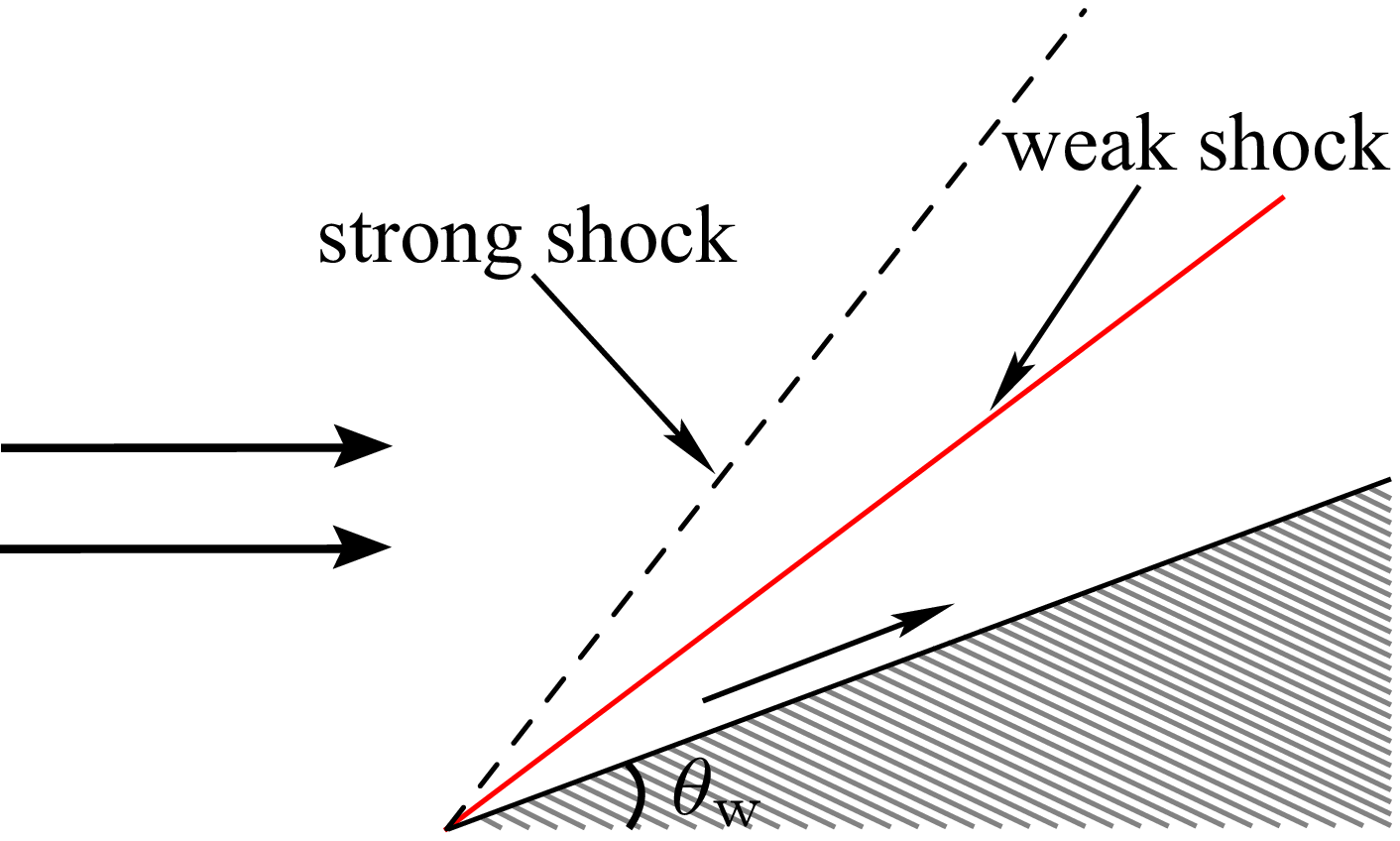}
		\caption{The shock polar in the $\uu$-plane and uniform steady (weak and strong) shock flows}
		\label{Figure2}
 \end{figure}

The fundamental issue -- whether one or both of the steady weak and strong shocks are physically
admissible -- has been vigorously debated over the past eight decades
({\it cf.} \cite{CourantF,Dafermos,Liu,Neumann,Serre}).
On the basis of experimental and numerical evidence, it has strongly indicated
that the steady weak shock
solution would be physically admissible, as Prandtl conjectured in \cite{Prandtl}.
One natural approach to single out the physically admissible steady
shock solutions is via the
stability analysis: The stable ones are physical.
For example, it is indicated in Courant-Friedrichs \cite{CourantF}, Section 123, that
``{\it The question arises which of the two actually occurs. It has frequently been stated that
the strong one is unstable and that, therefore, only the weak one could occur.
A convincing proof of this instability has apparently never been given.}''
On August 17, 1949, during the Symposium on the Motion of Gaseous Masses of Cosmical Dimensions held
at Paris, von Neumann \cite{Neumann} invited several eminent scientists of that days, including
Burges, Heisenberg, Liepmann, von K\'{a}rm\'{a}n, and Temple, to join his discussion panel on the topic entitled
``{\it the Existence and Uniqueness of Multiplicity of Solutions of the Aerodynamical Equations}''.
In his open remarks, von Neumann made his comments specifically on the wedge problem:
``{\it Occasionally the simplest hydrodynamical problems have several solutions, some of which are
very difficult to exclude on mathematical grounds only. For instance, a very simple hydrodynamical
problem is that of the supersonic flow of a gas through a concave corner, which obviously
leads to the appearance of shock wave. In general, there are two different solutions
with shock waves, and it is perfectly well known from experimentation that only one of the two,
the weaker shock wave, occurs in nature. But I think that all stability arguments
to prove that it must be so, are of very dubious quality.}''
From these comments, one may see that von Neumann was not so optimistic then
whether the mathematical stability analysis could provide a complete understanding of the non-uniqueness issue.
After von Neumannn's remarks, the panel members made their deep insights and different points of view
for the admissibility of
the weak and strong shock waves and related problems.
For more details, we refer the reader to
von Neumann \cite{Neumann}; see also \cite{Liu,Serre}.

It is interesting to observe that there had not too much progress on the global stability analysis of
the steady oblique shock solutions until recently; this is partly owing to the lack of
mathematical tools and  techniques that are required for solving the problem.
Mathematically, there are two levels of the stability analysis
of the multidimensional shock waves:
One is the static stability of the shocks under steady perturbations
within the steady regime,
{\it i.e.},
the steady perturbation of both the upstream supersonic flow and the wedge boundary;
the other is the dynamic stability to show that the steady shock solutions
are the long-time asymptotic limiting states of the corresponding unsteady solutions of the
Euler equations for compressible fluids.

As far as we have known, the rigorous study of the local static stability of supersonic shock waves
({\it i.e.}, both the upstream and downstream states are supersonic)
around the wedge vertex
for potential flow was first initiated by the Fudan Nonlinear PDE Group led by
Chaohao Gu and Tatsien Li in 1960; see \cite{GuLi}.
In this work, the shock wave involved was first regarded as a free boundary
to formulate the stability problem as a free boundary problem,
and was further reformulated the free
boundary problem into a fixed boundary problem via the hodograph transform.
The hodograph method was further refined by employing the parameterization method of
characteristics in  Gu-Li-Hou \cite{GuLiHou1,GuLiHou2,GuLiHou3}.
In \cite{Gu},
Gu extended the analysis of the wedge problem from the isentropic
to full Euler equations and related quasilinear hyperbolic systems.
In Li-Yu \cite{LiYu1,LiYu2} (see also \cite{LiYu3}),
a more general method of transformations
was developed to analyze more general quasilinear
hyperbolic systems.
The local stability of the
two-dimensional supersonic
shock waves for the full Euler equations was proved in Li \cite{Li-0,Li}
via employing the approach to solving free
boundary problems for quasilinear hyperbolic systems,
developed in Li-Yu \cite{LiYu1,LiYu2} (see also \cite{LiYu3}),
which considerably simplified
the original proof of  Schaeffer \cite{Schaeffer}
that had been achieved via the Nash-Moser iterations.
It is shown  that the flow possesses the
same qualitative features as a flow past a straight-sided wedge locally.
The result in Li \cite{Li-0,Li} was obtained
without  the  additional hypothesis  of the H\"{o}lder  continuity
employed in \cite{Schaeffer} via the Nash-Moser iterations.
See also Chen \cite{Chen0,Chen0-1}
for the local static stability of supersonic shock waves past a three-dimensional wing and
conical body for potential flow.
For the first rigorous treatment of
the  local existence  and  stability  of unsteady multidimensional
shock fronts for nonlinear hyperbolic systems of conservation laws,
see Majda \cite{Majda1,Majda2,Majda3}

The global stability results we present here are originally motivated by these
fundamental results, insights, and remarks mentioned above.
The purpose of this paper is to analyze the stability of both weak and strong
shock waves,
to show how the wedge problem can be formulated as mathematical
problems -- initial-boundary value type problems and free boundary problems,
to present some recent developments,
and to discuss further mathematical challenges and open problems
in the global stability analysis of multidimensional shock waves.

More precisely,
in Section 2, we first formulate the wedge problem as
an initial-boundary value type problem, and then reformulate it
into a free boundary problem.
In Section 3, we present the global static stability of steady supersonic
shocks under the BV perturbation of both the upstream steady supersonic flow
and the slope of the wedge boundary,
as long as the wedge vertex angle is less than the sonic angle.
In Section 4, we present the global static stability of both weak and strong transonic shocks
under the perturbation of both the upstream flow and the slope of the wedge boundary in
a weighted H\"{o}lder space.
In Section 5, we show that the steady weak supersonic/transonic shock solutions
are the asymptotic limits
of the dynamic self-similar solutions,
the Prandtl-Meyer configurations.
Finally, in Section 6, we discuss some further developments and mathematical challenges
in this research direction.

\section{Static Stability I: Mathematical Formulations and Free Boundary Problems}

In this section, we first formulate the wedge problem as
an initial-boundary value type problem, and then reformulate it
into a free boundary problem
when the perturbation of both the upstream steady supersonic flow and the
wedge boundary are suitably regular and small.

In order that
a piecewise smooth solution $U=(\uu, p, \rho)$
separated by a front $\mathcal{S}:=\{\xx\,:\,  x_2=\sigma(x_1), x_1\ge 0\}$
becomes a weak solution of the Euler equations \eqref{Euler1}, the
Rankine-Hugoniot conditions must be satisfied along
$\mathcal{S}$:
\begin{equation} \label{con-RH}
\begin{cases}
\sigma'(x_1)[\,\rho u_1\,]=[\,\rho u_2\,],\\[1mm]
\sigma'(x_1)[\,\rho u_1^2 + p\,]= [\,\rho u_1 u_2\,],\\[1mm]
\sigma'(x_1)[\,\rho u_1 u_2\,] =  [\,\rho u_2^2 + p\,],\\[1mm]
\sigma'(x_1)[\,\rho u_1(E+\frac{p}{\rho})\,] =[\,\rho u_2(E+\frac{p}{\rho})\,],
\end{cases}
\end{equation}
where $[\, \cdot\, ]$ denotes
the jump between the quantity of two states across front $\mathcal{S}$;
that is, if  $w^-$ and $w^+$ represent the left and right states, respectively,
then $[w]:= w^+ - w^-$.

Such a front $\mathcal{S}$ is called a shock if the entropy condition holds along $\mathcal{S}$:
{\it The density increases in the fluid direction across $\mathcal{S}$.}

\smallskip
For given a state $U^-$, all states $U$ that can be connected
with $U^-$ through the relations in
\eqref{con-RH}
form a curve in the state space $\R^4$;
the part of the curve whose states satisfy the entropy condition is called the shock polar.
The projection of the shock polar onto the $\uu$--plane is shown in Fig.~\ref{Figure2}.

In particular,  for an upstream uniform horizontal flow
$U_0^-=(u_{10}^-,0, p_0^-, \rho_0^-)$ past the upper part of a straight-sided wedge whose angle
is $\theta_{\rm w}$, the downstream constant flow can be
determined by
the shock
polar (see Fig.~\ref{Figure2}).
According to the shock polar, the
two flow angles are important: One is the detachment
angle $\theta_{\rm w}^{\rm d}$
that ensures the existence of an attached shock at the wedge
vertex, and the other is the sonic angle $\theta_{\rm w}^{\rm s}<\theta_{\rm w}^{\rm d}$
for which the
downstream fluid velocity is at the sonic speed in the direction.
More precisely, in Fig. \ref{Figure2},
$\theta_{\rm w}^{\rm s}$ is the wedge-angle such that
line $u_2=u_1 \tan\theta_{\rm w}^{\rm s}$ intersects
with the shock polar at a point on the circle of radius $c_0$,
and $\theta_{\rm w}^{\rm d}$ is the wedge-angle so that
line $u_2=u_1 \tan\theta_{\rm w}^{\rm d}$
is tangential to the shock polar and there is no intersection between
line $u_2=u_1 \tan \theta_{\rm w}$ and
the shock polar when $\theta_{\rm w}>\theta_{\rm w}^{\rm d}$.

When the wedge angle is less than
the detachment angle $\theta_{\rm w}^{\rm d}$,
the tangent point $T$ corresponding to the detachment angle
divides arc $\wideparen{HS}$ into the two open arcs
$\wideparen{TS}$ and $\wideparen{TH}$; see Fig. 1.2.
The nature of these two cases, as well as the case for arc $\wideparen{SQ}$,
is very different.
When the wedge angle $\theta_{\rm w}$ is between
$\theta_{\rm w}^{\rm s}$ and $\theta_{\rm w}^{\rm d}$,
there are two subsonic solutions; while the wedge angle
$\theta_{\rm w}$ is smaller than $\theta_{\rm w}^{\rm s}$,
there are one subsonic
solution and one supersonic solution.
Such an oblique shock $\mathcal{S}_0$ is also straight,
described by $x_2= s_0 x_1$.
The question is whether the steady oblique shock
solution is stable under a perturbation of both the upstream supersonic
flow and the wedge boundary.

Assume that the perturbed upstream flow $U^-_I$ is close to $U_0^-$,
which is supersonic and almost horizontal, and the wedge is close
to a straight-sided wedge. Then, for any suitable wedge angle (smaller
than a detachment angle), it is expected that there should be a shock
attached to the wedge vertex.
We now use a function $b(x_1)\ge 0$ to describe the upper wedge boundary:
\begin{equation}\label{wall}
\partial\mathcal{W}=\{\xx\in \R^2\,:\,  x_2= b(x_1),\  b(0)=0\}.
\end{equation}
Then the wedge problem can be formulated as the following problem:
\begin{problem}[Initial-Boundary Value Type Problem]
Find a global solution of system \eqref{Euler1}
in $\Omega:=\{x_2>b(x_1), x_1>0\}$
such that the following holds{\rm :}
\begin{enumerate}
\item[\rm (i)] Cauchy condition at $x_1=0${\rm :}
\begin{equation}
U|_{x_1=0}=U^-_I(x_2){\rm ;}
\end{equation}

\item[\rm (ii)] Boundary condition on $\partial\mathcal{W}$ as the slip boundary{\rm :}
\begin{equation}\label{slipcon}
\uu\cdot \mathbf{n}|_{\partial\mathcal{W}}=0,
\end{equation}
where $\mathbf{n}$ is the outer unit normal vector to $\partial \mathcal{W}$.
\end{enumerate}
\end{problem}

Suppose that the background shock is the straight line
given by $x_2= \s_0(x_1)=s_0 x_1$.
When the upstream steady supersonic
perturbation $U^-_I(x_2)$ at $x_1=0$
is suitably regular and small,
the upstream steady supersonic smooth solution $U^-(\xx)$ exists
in region $\Om^- = \Big\{\xx\,:\, 0<x_1< \frac{s_0}{2} x_1\Big\}$,
beyond the background shock, but is still close to $U_0^-$.

Suppose that the shock wave $\mathcal{S}$ we seek is
\[
\{\xx \,:\, \s(0)=0, \, x_2 = \s(x_1), x_1 \ge 0\}.
\]
The domain for the downstream flow behind $\mathcal{S}$ is denoted by
\begin{equation}\label{domain:1}
\Omega=\{\xx\in \R^2\,:\,  b(x_1)<x_2<\s(x_1), x_1>0\}.
\end{equation}

\begin{figure}
 \centering
\includegraphics[height=47mm]{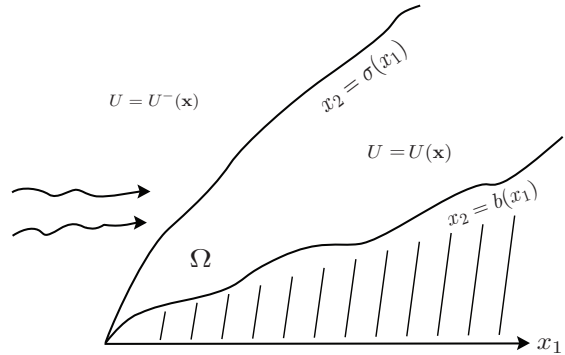}
\caption[]{The leading steady shock as a free boundary under the perturbation}
\end{figure}

Then {\bf Problem 2.1} can be further reformulated into
the following free boundary problem:

\begin{problem}[Free Boundary Problem; see also Fig. {\rm 2.1}]
Let $(U^-_0, U^+_0)$ be a constant transonic solution with transonic shock
$\mathcal{S}_0:=\{x_2=\sigma_0(x_1)=s_0x_1\}$.
For any upstream flow $U^-$ for equations \eqref{Euler1} in domain $\Om^-$
as a small perturbation of $U^-_0$,
find a shock $\CS:=\{x_2= \s(x_1)\}$ and a solution $U$ in $\Omega$ {\rm (}see Fig. {\rm 2.1}{\rm )},
which are small perturbations of $\mathcal{S}_0$ and $U^+_0$, respectively, such that
\begin{enumerate}
\item[\rm (i)] $U$ satisfies the equations in \eqref{Euler1} in domain $\Om${\rm ;}
\item[\rm (ii)]  The slip condition \eqref{slipcon} holds  along the wedge boundary $\del \mathcal{W}${\rm ;}
\item[\rm (ii)] The Rankine-Hugoniot conditions \eqref{con-RH} as free boundary conditions
hold along the shock front  $\CS$.
\end{enumerate}
When $U^+_0$ corresponding to a state
on arc $\wideparen{SQ}$  gives a weak supersonic
shock {\rm (}{\it i.e.}, both the upstream and downstream states are supersonic{\rm )} {\rm (}see Fig. {\rm 1.2}{\rm )},
the problem is denoted by {\rm \textbf{Problem 2.2(SS)}};
when $U^+_0$ corresponding to a subsonic state
on arc $\wideparen{TS}$ gives a weak transonic
shock {\rm (}{\it i.e.}, the upstream state is supersonic and the downstream
state is subsonic{\rm )} {\rm (}see Fig. {\rm 1.2}{\rm )},
the problem is denoted by {\rm \textbf{Problem 2.2 (WT)}};
while
the strong transonic shock problem corresponds to arc $\wideparen{TH}$,
denoted by  {\rm \textbf{Problem 2.2 (ST)}}.
\end{problem}

In general, the initial-boundary value type  problem ({\bf Problem 2.1}) is more general
than the free boundary problem ({\bf Problem 2.2}).
On the other hand, the complete solution
to the free boundary problem ({\bf Problem 2.2})
provides the global structural stability of
the steady oblique shocks, as well as more detailed structure of solutions.

\section{Static Stability II: Steady Supersonic Shocks}

If the downstream flow is supersonic ({\it i.e.},  $U_0^+\in \wideparen{SQ}$),
the corresponding shock is a weaker supersonic shock.

As indicated in \S 2, the rigorous study of the local static stability
of such shock waves
around the wedge vertex  for
the potential flow equation was first initiated by
the Fudan Nonlinear PDE Group led by Gu and Li in \cite{GuLi};
also see \cite{GuLiHou1,GuLiHou2,GuLiHou3,Gu,LiYu1,LiYu2,LiYu3}.
For the full Euler equations, the local stability of the supersonic shocks
was established by Gu \cite{Gu}, Li \cite{Li-0,Li}, and  Schaeffer \cite{Schaeffer} via different
approaches.

Global potential solutions were
constructed in \cite{Chen1,Chen2,Chen3,CXY,ChenFang,Zh1,Zh2} when the
wedge has certain convexity or the wedge is a small perturbation of
the straight-sided wedge with fast decay in the flow direction, whose
vertex angle is less than the detachment angle. In particular, in
Zhang \cite{Zh2}, the existence of two-dimensional steady supersonic
potential flows past piecewise smooth curved wedges, which are a
small perturbation of the straight-sided wedge, was established.

For the free boundary problem \textbf{Problem 2.2 (SS)},
the more general initial-boundary value type problem ({\bf Problem 2.1}) has been
solved for more general perturbations of both the initial data and wedge boundary,
in Chen-Zhang-Zhu  \cite{ChenZhangZhu} and Chen-Li \cite{ChenLi}.
More precisely,
\begin{enumerate}\renewcommand{\theenumi}{(\roman{enumi})}
\item The wedge boundary function $x_2=b(x_1)$ is  a Lipschitz function,
$b\in Lip(\R_+)$,
with
$$
TV(b'(\cdot))<\infty, \quad b'(0+)=0, \quad b(0)=0, \quad \arctan(b'(x_1))<\theta_{\rm w}^{\rm s},
$$
so that
${\bf n}(x_1\pm) =\frac{(-b'(x_1\pm),1)}{\sqrt{(b'(x_1\pm))^2+1}}$
is the outer normal vector to $\Gamma$ at point $x\pm$
(see Fig. 3.1);

\item
The upstream flow $U^-_I=(\uu^-_I,p^-_I,\rho^-_I)$  is a BV function ({\it i.e.}, $U^-_I\in BV(\R_+)$) satisfying that
$$
u^-_I>0,
\qquad |\uu^-_I|^2 > (c^-_I)^2:=\frac{\gamma p^-_I}{\rho^-_I}.
$$
\end{enumerate}

\begin{figure}
 \centering
\includegraphics[height=55mm]{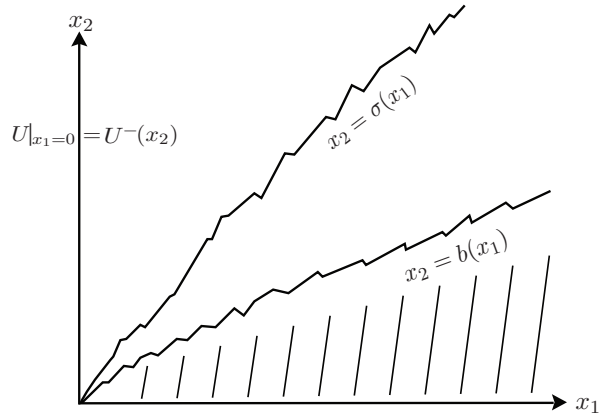}
\caption[]{The leading supersonic shock as a Lipschitz free boundary under the BV perturbations}
\end{figure}

With this setup, we have

\medskip
\noindent
\begin{theorem}[Chen-Zhang-Zhu \cite{ChenZhangZhu} and Chen-Li \cite{ChenLi}; see also Fig. 3.1]\quad
There are  $\e>0$ and $C>0$ such that, if
\begin{equation}\label{boundaryassump:1}
TV(U^-_I(\cdot))+TV (b'(\cdot))<\varepsilon,
\end{equation}
then there exists a pair of functions $(U(\xx), \sigma(x_1))${\rm :}
$$
U\in BV_{loc}(\R^2_+),\qquad \sigma'\in BV(\R_+),
$$
such that

{\rm (i)} The curve, $x_2=\sigma(x_1)$, is a leading shock above the wedge boundary $x_2=b(x_1)$ for any $x_1>0${\rm ;}

\smallskip
{\rm (ii)} $U$ is a global entropy solution of {\bf Problem 2.1}
in $\Omega:=\{\xx\,:\, x_2>b_1(x_1), x_1\ge 0\}$ with
\begin{eqnarray}
&&TV\{ U(x_1,\cdot)\,:\, (b(x_1), \sigma(x_1))\cup (\sigma(x_1), \infty)\} \le C\big(TV (U^-_I(\cdot))+ TV (b'(\cdot))\big)
\,\,\,\,\, \mbox{for any}\, \, x_1 \in \R_+,\qquad \label{1.11a}\\[1mm]
&& (u, v)\cdot {\bf n}|_{x_2=b(x_1)}=0
\quad\mbox{in the trace sense}; \label{1.12a}
\end{eqnarray}

{\rm (iii)} There exist constants $s_{\infty}$ and $p_{\infty}$ such that
$$
\lim_{x_1\rightarrow\infty}|\sigma'(x_1)-s_{\infty}| =0, \quad \lim_{x_1\rightarrow\infty}
\sup\big\{|p(\xx)-p_{\infty}|\, :\; b(x_1)<x_2<\sigma(x_1)\big\} =0,
$$
and
$$
\lim_{x_1\rightarrow\infty}
\sup\big\{\big|\frac{u_2(\xx)}{u_1(\xx)}
  -b'(\infty)\big|\,: \; b(x_1)<x_2< \sigma(x_1)\big\} =0.
$$

Moreover, the entropy solution $U=U(\xx)$ is stable with respect to the initial $BV$ perturbation
in $L^1$ and unique in a broader class -- the class of viscosity solutions.
\end{theorem}

This theorem indicates that the leading steady supersonic oblique shock-front emanating from
the wedge vertex is nonlinearly stable in structure, although there
may be many weaker waves and vortex sheets between the leading supersonic shock-front
and the wedge boundary or the $x_2$--axis where the initial condition is assigned,
under the $BV$ perturbation of both the upstream flow
and the slope of the wedge boundary, as long as the wedge vertex angle is
less than the sonic angle $\theta_{\rm w}^{\rm s}$.
Moreover, the steady supersonic shock for the wedge problem
is nonlinearly stable in $L^1$ under the $BV$ perturbation.
This asserts that
the steady supersonic oblique shock should be
physical admissible, as observed from the experimental results.

More specifically,
in Chen-Zhang-Zhu \cite{ChenZhangZhu}, in order to establish the global existence
of solutions for the constant Cauchy data $U^-_I=U_0^-$, we first developed
a modified Glimm
scheme and identified a Glimm-type functional by incorporating the
curved wedge boundary and the strong shock naturally, and by
tracing the interactions not only between the wedge boundary and
weak waves but also the interaction between the strong shock
and weak waves. Some detailed interaction estimates are carefully
made to ensure that the Glimm-type functional monotonically
decreases in the flow direction. In particular, one of the essential
estimates is on the strengths of the reflected waves for
system \eqref{Euler1} in the interaction
between the strong shock and weak waves; and the second
essential estimate is the interaction estimate between the wedge
boundary and weak waves. Another essential estimate is on tracing
the approximate strong shocks in order to establish the nonlinear
stability and asymptotic behavior of the strong shock
emanating from the wedge vertex under the wedge perturbation.

In Chen-Li \cite{ChenLi}, based on the understanding
of the problem in Chen-Zhang-Zhu \cite{ChenZhangZhu},
we further established
the $L^1$ well-posedness for \textbf{Problem 2.1} when
the total variation of
both the boundary slope function and
the Cauchy data (upstream flow) is small.
We
first obtained the existence of solutions in $BV$ when the upstream
flow $U^-_I$ has small total variation by the wave front tracking method and
then established the $L^1$--stability of the solutions with respect to
the upstream flows. In particular, we incorporated the nonlinear
waves generated from the wedge boundary to develop a Lyapunov
functional between two solutions containing the strong shock fronts,
which is equivalent to the $L^1$--norm, and proved that the functional
decreases in the flow direction. Then the $L^1$--stability was
established, which implies the uniqueness of the solutions by the wave front
tracking method. Finally, the uniqueness of solutions in a broader
class, the class of viscosity solutions, was also obtained.

\section{Static Stability III: Weak and Strong Transonic Shocks}

For transonic ({\it i.e.}, supersonic-subsonic) shocks,
there are two cases -- the transonic shock with the subsonic state corresponding
to arc $\wideparen{TS}$ (which is a weaker shock) and the one corresponding
to arc $\wideparen{TH}$ (which is a stronger shock) (see Fig.~1.2).
The strong shock case was first studied in Chen-Fang \cite{ChenFang} for the potential flow.
In Fang \cite{Fang},  the full Euler equations were studied
with a uniform Bernoulli constant
for both weak and strong transonic shocks.
Because the framework is a weighted Sobolev space, the asymptotic behavior of the shock
slope or subsonic solution was not derived.
In Yin-Zhou \cite{YinZhou}, the H\"older norms were used for
the estimates of solutions of the full Euler equations
with the assumption on the sharpness of the wedge angle,
which means that the subsonic state
is near point $H$ in the shock polar.
In Chen-Chen-Feldman \cite{CCF2}, the weaker transonic shock, which corresponds to
arc $\wideparen{TS}$, was investigated; and
the existence, uniqueness, stability, and asymptotic
behavior of  subsonic solutions were obtained.
In \cite{CCF2,YinZhou}, a potential function is used to reduce the full Euler equations to
one elliptic equation in the subsonic region.
The method was first proposed in \cite{CCF} and
has the advantage of integrating the conservation properties of the Euler system
into a single elliptic
equation. However, working on the potential function requires at least the Lipschitz estimate
of the potential function to keep the subsonicity of the flow.
In our recent paper \cite{CCF3}, we have directly employed the decomposition of the full Euler equations
into two algebraic equations
and a first-order elliptic system of two equations and have established the stability and asymptotic
behavior of transonic flows for {\bf Problem 2.2 (WT)--(ST)} in a weighted H\"{o}lder space.

To state the results, following \cite{CCF3},  we need to introduce the weighed H\"older norms
in the subsonic domain $\Omega$,
where $\Omega$ is either a truncated triangular domain or an unbounded domain with the vertex at
origin $\O$ and one side as the wedge boundary.
There are two weights: One is the distance function to origin $\O$ and the other is to
the wedge boundary $\partial \mathcal{W}$.
For any $\xx, \xx'\in \Omega$,  define
\begin{align*}
&\gd^{\rm o}_\xx := \min(|\xx|,1),\quad \gd^{\rm o}_{\xx,\xx'} := \min (\gd^{\rm o}_\xx, \gd^{\rm o}_{\xx'}),
\quad \gd^{\rm w}_\xx := \min(\textrm{dist}(\xx,\partial \mathcal{W}),1), \quad \gd^{\rm w}_{\xx,\xx'}
   := \min (\gd_\xx^{\rm w},\gd_{\xx'}^{\rm w}), \\[1mm]
&\gD_\xx := |\xx|+1, \quad \gD_{\xx, \xx'}:=\min(\gD_\xx, \gD_{\xx'}),
\quad \widetilde{\gD}_\xx := \textrm{dist}(\xx,\partial \mathcal{W})+1,
\quad \widetilde{\gD}_{\xx, \xx'}:=
\min(\widetilde{\gD}_\xx, \widetilde{\gD}_\xx).
\end{align*}
Let $\ga \in (0,1)$,  $\tau,l, \gamma_1, \g_2 \in \R$ with $\g_1 \ge \g_2$,
and let $k$ be a nonnegative integer.
Let $\kk = (k_1, k_2)$ be an integer-valued vector, where $k_1, k_2 \ge 0$,
$|\kk|=k_1 +k_2$, and $D^{\kk}= \po_{x_1}^{k_1}\po_{x_2}^{k_2}$.
We define
\begin{align}
&[ f ]_{k,0;(\tau,l); \Omega}^{(\g_1;\O)(\g_2;\partial \mathcal{W})}
=\sup_{\begin{subarray}{c}
\xx\in \Omega\\
 |\kk|=k
\end{subarray}}
\begin{array}{l}
\big\{(\gd^{\rm o}_{\xx})^{\max\{\g_1 + \min\{k,-\g_2\},0\}}(\gd^{\rm w}_{\xx})^{\max\{k+\g_2,0\}}
\,\gD_\xx^{\tau} \widetilde{\gD}_\xx^{l+k} |D^\kk f(\xx)|\big\},
\end{array}\label{def-normk0}\\[2mm]
&{[ f ]}_{k,\ga;(\tau,l); \Omega}^{(\g_1;\O)(\g_2;\partial \mathcal{W})}
\nonumber\\&\quad
=
 \sup_{
\begin{subarray}{c}
   \xx, \xx'\in \Omega\\
 \xx \ne \xx', |\kk|=k
 \end{subarray}}
\Big\{
  (\gd^{\rm o}_{\xx,\xx'})^{\max\{\g_1 + \min\{k+\ga,-\g_2\},0\}}(\gd^{\rm w}_{\xx,\xx'})^{\max\{k+\ga+\g_2,0\}}
  \gD_{\xx,\xx'}^{\tau} \widetilde{\gD}_{\xx,\xx'}^{l+k+\ga}\frac{|D^\kk f(\xx)-D^\kk
f(\xx')|}{|\xx-\xx'|^\ga}
\Big\},
\label{def-normkga}
\\
& \|f\|_{k,\ga;(\tau,l); \Omega}^{(\g_1;\O)(\g_2;\partial \mathcal{W})}= \sum_{i=0}^k {[ f
]}_{i,0 ;(\tau,l);\Omega}^{(\g_1;\O)(\g_2;\partial \mathcal{W})} + {[ f
]}_{k,\ga;(\tau,l);\Omega}^{(\g_1;\O)(\g_2;\partial \mathcal{W})}. \label{def-norm}
\end{align}

For a vector-valued function $\ff=(f_1, f_2, \cdots,f_n )$, we
define
\[
 \|\ff\|_{k,\ga;(\tau,l); \Omega}^{(\g_1;\O)(\g_2;\partial \mathcal{W})} =
\sum_{i=1}^n  \|f_i\|_{k,\ga;(\tau,l); \Omega}^{(\g_1;\O)(\g_2;\partial \mathcal{W})}.
\]
 Let
\begin{equation}\label{def-C}
C^{k,\ga;(\tau,l)}_{(\g_1;\O)(\g_2;\partial \mathcal{W})}(\Omega)
= \left\{ f: \|f\|_{k,\ga;(\tau,l); \Omega}^{(\g_1;\O)(\g_2;\partial \mathcal{W})} <\infty \right\}.
\end{equation}

The requirement that $\g_1 \ge \g_2$ in the definition above means that the regularity up
to the wedge boundary is no worse than the regularity up to the wedge vertex.
When $\g_1 = \g_2$, the $\d^{\rm o}$--terms disappear so that $(\g_1,\O)$ can be dropped
in the superscript.
If there is no weight $(\g_2,\pw)$ in the superscript, the $\d$--terms for the weights
should be understood as $(\d_{\xx}^{\rm o})^{\max\{k+\g_1,0\}}$
and $(\d_{\xx}^{\rm o})^{\max\{k + \ga +\g_1,0\}}$
in \eqref{def-normk0} and \eqref{def-normkga}, respectively.
Moreover, when no weight appears in the superscripts of the seminorms
in \eqref{def-normk0} and \eqref{def-normkga},
it means that neither $\delta^{\rm o}$ nor $\delta^{\rm w}$ is present.
For a function of one variable defined on $(0,\infty)$,
the weighted norm $\|f\|^{(\g_2;0)}_{k,\ga;(l);\R^+}$
is understood in the same as the definition
above with the weight to $\{0\}$ and the decay
at infinity.

Since the variables in $U$
are expected to have different levels of regularity,
we distinguish these variables by defining
${U_1}=(u_1,\rho)$ and $U_2=(w, p)$ for $w=\frac{u_2}{u_1}$.
Let $U_{10}^+$ and $U_{20}^+$ be
the corresponding background subsonic states.

\begin{theorem}[Chen-Chen-Feldman \cite{CCF3}] \label{thm-main} \quad
There are positive constants $\ga , \b,  C_0$, and $\varepsilon$, depending only on the background
states $(U^-_0, U^+_0)$,  such that there exists a solution $(U,\s)$
for either of {\rm \textbf{Problem 2.2 (WT)}} and {\rm \textbf{Problem 2.2 (ST)}}
such that $(U,\s)$ satisfies the following
estimates{\rm :}

\begin{enumerate}
\item[\rm (i)] For {\rm \textbf{Problem 2.2 (WT)}},
\begin{equation}\label{est-U-small-pert}
\begin{split}
& \|U_1 - U_{10}^+\|^{(-\ga;\del  \mathcal{W})}_{2,\ga;(0,1+\gb);\Om}
+\|U_2 - U_{20}^+\|^{(-\ga;  \O)(-1-\ga;\del  \mathcal{W})}_{2,\ga;(1+\gb,0);\Om}
+ \| \s' - s_0\|^{(-\ga;0)}_{2,\ga;(1+\gb);\R^+}\\
&\le{}  C_0 \left(\| U^- -U^-_0 \|_{2,\ga;(1+\b,0);\Om^-}
     +\|b'\|^{(-\ga;0)}_{1,\ga; (1+\b);\R^+} \right),
\end{split}
\end{equation}
provided that
\[
\| U^- -U^-_0 \|_{2,\ga;(1+\b,0);\Om^-} +\| b' \|^{(-\ga;0)}_{1,\ga; (1+\b);\R^+} <\ve;
\]

\item[\rm (ii)] For {\rm \textbf{Problem 2.2 (ST)}},
\begin{equation}\label{est-U-small-pert2}
\begin{split}
&\|U_1 - U_{10}^+\|^{(-1-\ga;\pw)}_{2,\ga;(0,\gb);\Om}
	+\|U_2 - U_{20}^+\|^{(-1-\ga; \O)}_{2,\ga;(\gb,0);\Om}
+\|\s' - s_0\|^{(-1-\ga;0)}_{2,\ga;(\gb);\R^+}
\\
&\le{} C_0  \left(\| U^- -U^-_0\|_{2,\ga;(\b);\Om^-}
  +\|b' \|^{(-1-\ga;0)}_{2,\ga; (\b);\R^+}  \right),
\end{split}
\end{equation}
provided that
\[
\|U^- -U^-_0\|_{2,\ga;(\b,0);\Om^-} + \|b'\|^{(-\ga -1;0)}_{2,\ga; (\b);\R^+}<\ve.
\]
\end{enumerate}
The solution $(U,\s)$ is unique within the class of solutions such that
  the left-hand side of
\eqref{est-U-small-pert} for {\rm \textbf{Problem 2.2 (WT)}} or \eqref{est-U-small-pert2}
for {\rm \textbf{Problem 2.2 (ST)}} is less than $C_0 \vae$.
\end{theorem}

The dependence of constants $\ga , \b,  C_0$, and $\varepsilon$ in Theorem 4.1
is as follows:  $\ga$ and $\b$ depend on $U^-_0$ and $U^+_0$,
but are independent of $C_0$ and $\varepsilon$;
$C_0$ depends on $U^-_0, U^+_0,\ga$, and $b$, but are independent of $\ve$;
and $\ve$ depends on all $U^-_0, U^+_0, \ga , \b$, and $C_0$.

The difference in the results of the two problems is that the solution of
{\rm \textbf{Problem 2.2 (WT)}} has less regularity at corner ${O}$
and decays faster with respect to $|\xx|$ (or the distance from the wedge boundary)
than the solution of {\rm \textbf{Problem 2.2 (ST)}}.

To achieve these, the main strategy is to use the physical variables
to do the estimates, instead of the potential function.
The advantage of this method is that only the lower regularity ({\it i.e.}, the $C^0$--estimate)
is enough to guarantee the subsonicity.
Furthermore, directly estimating the physical state function $U$
also yields a better asymptotic decay rate:
For weaker transonic shocks, in our earlier paper \cite{CCF2},
the decay rate of the velocity is only $|\xx|^{-\gb}$;
while, as indicated in \eqref{est-U-small-pert} here,
the subsonic solution decays to a limit state at
rate $|\xx|^{-\gb-1}$.

More precisely, we have first used the Lagrangian coordinates to straighten the streamlines.
The reason for this is that the Bernoulli variable and entropy are conserved
along the streamlines.
Using the streamline as one of the coordinates simplifies the formulation,
especially for the asymptotic behavior of the solution.
Then
we have decomposed the Euler system into two algebraic equations and two elliptic equations, as in \cite{schen3,Fang}.
Differentiating the two elliptic equations yields a second-order elliptic equation
in divergence form for the flow direction $w= \frac{u_2}{u_1}$.
Given $U$ for the coefficients in the equation, we can solve for a new variable $\tilde{w}$.
Once we have solved for $\tilde{w}$ and obtained the desired estimates,
the rest variables have then been updated
so that a map $\d \widetilde{U}=\mathcal{Q}(\d U)$ has been constructed,
where $\d U$ and $\d \widetilde{U}$ are the perturbations from the background subsonic state.
The estimates based on our method do not yield the contraction for $\mathcal{Q}$.
Therefore, the Schauder fixed point argument has been employed to obtain the existence of the subsonic
solution. For the uniqueness, we have taken the difference of two solutions and have estimated the difference
by using the weighted H\"older norms with less decay rate.

One point we should emphasize here is that the decay pattern is different from the potential flow.
In a potential flow, the decay is with respect to $|\xx|$.
For example, if $\varphi$ converges to $\varphi_0$ at rate $|\xx|^{-\gb}$,
then $\grad \varphi$ converges at rate $|\xx|^{-\gb-1}$.
For the Euler equations, because the Bernoulli variable and the entropy are constant along
the streamlines,
$U_1=(u_1, \rho)$ does not converge to the background state along the streamlines,
but does converge only across the streamlines away from the wedge.
Therefore, when the elliptic estimates are performed,
the scaling is with respect to the distance from the wedge,
rather than  $|\xx|$.
This results in the following decay pattern:
In Lagrangian coordinates $\yy$, there exists an asymptotic limit $U^\infty = (u_1^\infty, 0, p_0^+, \rho^\infty)$;
$U$ converges to $U^\infty$ at rate $|\yy|^{-\gb}$, but $\grad U$ converges at rate $|\yy|^{-\gb} (y_2+1)^{-1}$.
That is, the extra decay for the gradient of the solution is only along the $y_2$--direction.

It would also be interesting to investigate the stability problems under more general perturbations, say the BV perturbation of
both the upstream flow and
the slope of the wedge boundary.

\section{Dynamic Stability: Self-Similar Transonic Shocks and Existence of Prandtl-Meyer Configurations
   for Potential Flow}

Since both weak and strong steady shock solutions are stable in the steady regime,
the static stability analysis alone
is not able to single out one of them in this sense,
unless an additional condition is posed on the speed of the downstream flow at infinity.
Then the dynamic stability analysis becomes more significant
to understand the non-uniqueness issue of the steady oblique shock solutions.
However, the problem for the dynamic stability of the steady shock solutions for
supersonic flow past solid wedges involves several additional
mathematical difficulties.
The recent efforts have been focused on the construction of the global Prandtl-Meyer
configurations in the self-similar coordinates for potential flow.

The compressible potential flow is governed by the conservation law of mass and the Bernoulli law:
\begin{align}
\label{1-a}
&\der_t\rho+ \nabla_{\bf x}\cdot (\rho \nabla_{\bf x}\Phi)=0,\\
\label{1-b}
&\der_t\Phi+\frac 12|\nabla_{\bf x}\Phi|^2+h(\rho)=B
\end{align}
for density $\rho$ and velocity-potential $\Phi$,
where $B$ is the Bernoulli constant determined by the upstream flow and/or boundary conditions,
and $h(\rho)$ is given by
$h(\rho)
=\int_1^{\rho}\frac{p'(\varrho)}{\varrho}\,d\varrho=\int_1^{\rho}\frac{c^2(\varrho)}{\varrho}\,d\varrho$
for the sound speed $c(\rho)$ and pressure $p$.
For an ideal polytropic gas, by scaling without loss of generality, the sound speed $c$ and pressure $p$ are given by
\begin{equation}
\label{cont-rel-pt}
  p(\rho)=\frac{\rho^{\gam}}{\gam}, \qquad   c^2(\rho)=\rho^{\gam-1},
  \qquad  h(\rho)=\frac{\rho^{\gam-1}-1}{\gam-1}
\end{equation}
for the adiabatic component $\gam>1$.

By \eqref{1-b}--\eqref{cont-rel-pt}, $\rho$ can be expressed as
\begin{equation}
\label{1-b1}
\rho(\der_t\Phi,\nabla_{\bf x}\Phi)=h^{-1}\bigl(B-\der_t\Phi-\frac 12|\nabla_{\bf x}\Phi|^2\bigr).
\end{equation}
Then system \eqref{1-a}--\eqref{1-b} can be rewritten as
\begin{equation}
\label{1-b2}
\der_t\rho(\der_t\Phi, \nabla_{\bf x}\Phi)
+\nabla_{\bf x}\cdot\big(\rho(\der_t\Phi, \nabla_{\bf x}\Phi)\nabla_{\bf x}\Phi\big)=0
\end{equation}
with $\rho(\der_t\Phi, \nabla_{\bf x}\Phi)$ determined by \eqref{1-b1}.

As we discussed earlier, if a supersonic flow
with a constant density $\irho >0$ and a velocity ${\bf u}_{0}=(\iu, 0)$,
$\iu> c_0:=c(\rho_0)$,
impinges toward wedge $W$ in \eqref{wedge-1},
and if $\theta_{\rm w}$ is less
than the detachment angle $\theta_{\rm w}^{\rm d}$,
then the well-known {\emph{shock polar analysis}} shows that there are two different
steady weak solutions:
{\emph{the steady weak shock solution}} $\bar{\Phi}$ and {\emph{the steady strong shock solution}},
both of which satisfy the entropy condition and the slip boundary condition (see Fig. 5.3).

Then the dynamic stability of the weak transonic shock solution for potential flow can be
formulated
as the following problem:

\begin{problem}[Initial-Boundary Value Problem]
\label{problem-1}
Given $\gam>1$, fix $(\irho, \iu)$ with $\iu>c_0$.
For a fixed $\theta_{\rm w}\in (0,\theta_{\rm w}^{\rm d})$,
let $W$ be given by \eqref{wedge-1}.
Seek a global weak solution $\Phi\in W^{1,\infty}_{loc}(\R_+\times (\R^2\setminus W))$
of Eq. \eqref{1-b2} with $\rho$ determined by \eqref{1-b1}
and $B=\frac{\iu^2}{2}+h(\irho)$ so that $\Phi$ satisfies the initial condition at $t=0${\rm :}
\begin{equation}\label{1-d}
(\rho,\Phi)|_{t=0}=(\irho, \iu x_1) \qquad\,\,\, \text{for}\;\;\xx\in \R^2\setminus W,
\end{equation}
and the slip boundary condition along the wedge boundary $\der W${\rm :}
\begin{equation}\label{1-e}
\nabla_{\bf x}\Phi\cdot {\bf n}_{\rm w} |_{\der W}=0,
\end{equation}
where ${\bf n}_{\rm w}$ is the exterior unit normal to $\der W$.

In particular, we seek a solution $\Phi\in W^{1,\infty}_{loc}(\R_+\times (\R^2\setminus W))$
that converges to
the steady weak oblique shock solution $\bPhi$ corresponding to
the fixed parameters $(\irho, \iu, \gam, \theta_{\rm w})$
with $\bar{\rho}=h^{-1}(B-\frac 12|\nabla \bPhi|^2)$,
when $t\to \infty$, in the following sense{\rm :}
For any $R>0$, $\Phi$ satisfies
\begin{equation}
\label{time-asymp-lmt}
\lim_{t\to \infty} \|(\nabla_{\bf x}\Phi(t,\cdot)-\nabla_{\bf x}\bPhi,
\rho(t,\cdot)-\bar{\rho})\|_{L^1(B_R({\bf 0})\setminus W)}=0
\end{equation}
for $\rho(t,{\bf x})$ given by \eqref{1-b1}.
\end{problem}

\medskip
Since the initial data function in \eqref{1-d} does not satisfy the boundary condition \eqref{1-e},
a boundary layer is generated along the wedge boundary starting at $t=0$,
which forms the Prandtl-Meyer configuration, as proved in
Bae-Chen-Feldman \cite{BCF-14,BCF-16}.

\smallskip
Notice that the initial-boundary value problem, {\bf Problem 5.1},
is invariant under the scaling:
\begin{equation*}
(t, {\bf x})\rightarrow (\alp t, \alp{\bf x}),\quad (\rho, \Phi)\rightarrow (\rho, \frac{\Phi}{\alp})\qquad\quad \text{for}\;\;\alp\neq 0.
\end{equation*}
Thus, we seek self-similar solutions in the form of
\begin{equation*}
\rho(t, {\bf x})=\rho(\xxi),\quad \Phi(t, {\bf x})=t\phi(\xxi)\qquad\quad \text{for}\;\;\xxi=\frac{{\bf x}}{t}.
\end{equation*}
Then the pseudo-potential function $\vphi=\phi-\frac 12|\xxi|^2$ satisfies
the following Euler equations for self-similar solutions:
\begin{align}
\label{1-r}
{\rm div}(\rho D\vphi)+2\rho=0,\qquad
\frac{\rho^{\gam-1}-1}{\gam-1}+\big(\frac 12|D\vphi|^2+\vphi\big)=B,
\end{align}
where the divergence $div$ and gradient $D$ are with respect to $\xxi$.
From this, we obtain the following equation for the pseudo-potential function $\vphi(\xxi)$:
\begin{equation}
\label{2-1}
{\rm div}(\rho(|D\vphi|^2,\vphi)D\vphi)+2\rho(|D\vphi|^2,\vphi)=0
\end{equation}
for
\begin{equation}
\label{1-o}
\rho(|D\vphi|^2,\vphi)=
\bigl(B_0-(\gam-1)(\frac 12|D\vphi|^2+\vphi)\bigr)^{\frac{1}{\gam-1}},
\end{equation}
where we have set $B_0:=(\gam-1)B+1$.
Then we have
\begin{equation}
\label{1-a1}
c^2(|D\vphi|^2,\vphi)=
B_0-(\gam-1)\big(\frac 12|D\vphi|^2+\vphi\big).
\end{equation}
Equation \eqref{2-1} is an equation of mixed elliptic-hyperbolic type. It is elliptic if and only if
\begin{equation}
\label{1-f}
|D\vphi|<c(|D\vphi|^2,\vphi).
\end{equation}

As the upstream flow has the constant velocity $(\iu,0)$, the corresponding pseudo-potential
$\ivphi$ has the expression of
\begin{equation}
\label{1-m}
\ivphi=-\frac 12|\xxi|^2+\iu\xi_1.
\end{equation}
{\bf Problem 5.1} can then be reformulated as the following boundary value problem
in the domain:
$$
\Lambda:=\R^2_+\setminus\{\xxi\,: \,\xi_2\le \xi_1\tan\theta_{\rm w},\, \xi_1\ge 0\}
$$
in the
self-similar coordinates $\xxi$, which corresponds to domain $\{(t, {\bf x})\, :\, {\bf x}\in \R^2_+\setminus W,\, t>0\}$
in the $(t, {\bf x})$--coordinates.

\medskip
\begin{problem}[Boundary Value Problem]
\label{problem-5.2}
Seek a solution $\vphi$ of equation \eqref{2-1} in the self-similar domain $\Lambda$
with the slip boundary condition{\rm :}
\begin{equation}
\label{1-k}
D\vphi\cdot \mathbf{n}|_{\partial \Lambda}=0,
\end{equation}
and the asymptotic boundary condition at infinity{\rm :}
\begin{equation}\label{1-k-b}
\vphi-\vphi_0\longrightarrow 0
\end{equation}
along each ray $R_\theta:=\{ \xi_1=\xi_2 \cot \theta, \xi_2 > 0 \}$
with $\theta\in (\theta_{\rm w}, \pi)$ as $\xi_2\to \infty$
in the sense that
\begin{equation}\label{1-k-c}
\lim_{r\to \infty} \|\varphi  - \varphi_0\|_{C(R_{\theta}\setminus B_r(0))} = 0.
\end{equation}
\end{problem}

\begin{figure}[htp]
\centering
\begin{psfrags}
\psfrag{ls}[cc][][0.8][0]{$\leftshockop$}
\psfrag{sn}[cc][][0.8][0]{$\rightshockop$}
\psfrag{ol}[cc][][0.8][0]{$\oOmop$}
\psfrag{on}[cc][][0.8][0]{$\nOmop$}
\psfrag{lsn}[cc][][0.8][0]{$\phantom{aa}\leftsonicop$}
\psfrag{rsn}[cc][][0.8][0]{$\rightsonicop\phantom{aaaaaaa}$}
\psfrag{tw}[cc][][0.8][0]{$\theta_w$}
\psfrag{o}[cc][][0.8][0]{$O$}
\psfrag{lb}[cc][][0.8][0]{$\leftbottom$}
\psfrag{rb}[cc][][0.8][0]{$\rightbottom$}
\psfrag{lt}[cc][][0.8][0]{$\lefttop$}
\psfrag{rt}[cc][][0.8][0]{$\righttop$}
\psfrag{s}[cc][][0.8][0]{$\shock$}
\psfrag{om}[cc][][0.8][0]{$\Om$}
\includegraphics[scale=.5]{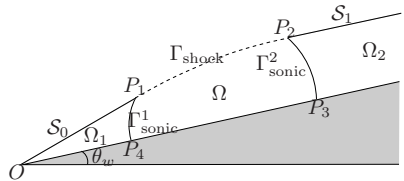}
\vspace{-24pt}
\caption{Self-similar solutions for $\theta_{\rm w}\in (0, \theta_{\rm w}^{\rm s})$
in the self-similar coordinates $\xxi$;
see Bae-Chen-Feldman \cite{BCF-14}}\label{fig:global-structure-1}
\end{psfrags}
\end{figure}

\begin{figure}[htp]
\centering
\begin{psfrags}
\psfrag{oi}[cc][][0.8][0]{$\phantom{aaaaaaa}\Om_{\infty}: (\iu, 0), \rho_\infty$}
\psfrag{om}[cc][][0.8][0]{$\Om$}
\psfrag{on}[cc][][0.8][0]{$\nOmop$}
\psfrag{sn}[cc][][0.8][0]{$\rightshockop$}
\psfrag{ls}[cc][][0.8][0]{$\leftshockop$}
\psfrag{s}[cc][][0.8][0]{\phantom{aaaaa}$\shock$}
\psfrag{tw}[cc][][0.8][0]{$\theta_w$}
\includegraphics[scale=0.7]{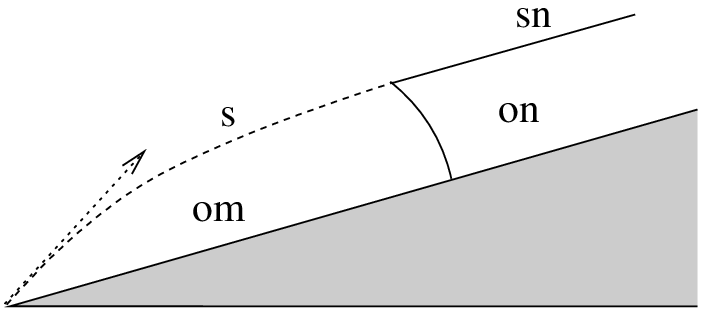}
\caption{Self-similar solutions for $\theta_{\rm w}\in [\theta_{\rm w}^{\rm s},\theta_{\rm w}^{\rm d})$
in the self-similar coordinates $\xxi$;
see Bae-Chen-Feldman \cite{BCF-14}}\label{fig:global-structure-2}
\end{psfrags}
\end{figure}

In particular, we seek a weak solution of {\bf Problem 5.2}
with two types of Prandtl-Meyer configurations whose occurrence
is determined by the wedge angle $\theta_{\rm w}$ for the two different cases:
One contains a straight weak oblique shock attached to the wedge vertex $O$, and
the oblique shock is connected to a normal shock through a curved shock
when $\theta_{\rm w}<\theta_{\rm w}^{\rm s}$,
as shown in Fig. \ref{fig:global-structure-1}; the other contains a curved shock attached
to the wedge vertex
and connected to a normal shock when $\theta_{\rm w}^{\rm s}\le \theta_{\rm w}<\theta_{\rm w}^{\rm d}$,
as shown in Fig. \ref{fig:global-structure-2}, in which
the curved shock $\Gamma_{\rm shock}$ is tangential to
a straight weak oblique shock $S_0$ at the wedge vertex.

\medskip

A shock is a curve across which $D\vphi$ is discontinuous. If $\Om^+$ and $\Om^-(:=\Om\setminus \ol{\Om^+})$
are two nonempty open subsets of $\Om\subset \R^2$, and $\mathcal{S}:=\der\Om^+\cap \Om$ is a $C^1$-curve
where $D\vphi$ has a jump, then $\vphi\in W^{1,1}_{loc}\cap C^1(\Om^{\pm}\cup S)\cap C^2(\Om^{\pm})$
is a global weak solution of \eqref{2-1} in $\Om$ if and only if $\vphi$ is in $W^{1,\infty}_{loc}(\Om)$
and satisfies equation \eqref{2-1} and the Rankine-Hugoniot condition on $\mathcal{S}$:
\begin{equation}
\label{1-h}
[\rho(|D\vphi|^2, \vphi)D\vphi\cdot\nnu]_{\mathcal{S}}=0,
\end{equation}
where $[F]_{\mathcal{S}}$ is defined by
$$
[F(\xxi)]_{\mathcal{S}}:=F(\xxi)|_{\overline{\Om^-}}-F(\xxi)|_{\overline{\Om^+}}\qquad\quad\text{for}\;\; \xxi\in \mathcal{S}.
$$
Note that the condition, $\vphi\in W^{1,\infty}_{loc}(\Om)$, requires that
\begin{equation}
\label{1-i}
[\vphi]_{\mathcal{S}}=0.
\end{equation}
The front, $\mathcal{S}$, is called a shock if density $\rho$ increases in the flow direction across $\mathcal{S}$.
A piecewise smooth solution whose discontinuities are all shocks is called an entropy solution.

To seek a global entropy solution of {\bf Problem \ref{problem-5.2}} with the structure
of Fig. \ref{fig:global-structure-1} or Fig. \ref{fig:global-structure-2},
one needs to compute the pseudo-potential function $\vphi_0$ below $S_0$.

\begin{figure}[htp]
\centering
\begin{psfrags}
\psfrag{tc}[cc][][0.8][0]{$\phantom{aaaaaaaaaaa}u_2=u_1\tan\theta_{\rm w}^{\rm s}$}
\psfrag{tw}[cc][][0.8][0]{$\phantom{aaaaaaaaaaa}u_2=u_1 \tan\theta_w$}
\psfrag{td}[cc][][0.8][0]{$\phantom{aaaaaaaaaaa}u_2=u_1\tan\theta_{\rm w}^{\rm d}$}
\psfrag{u}[cc][][0.8][0]{$u_1$}
\psfrag{v}[cc][][0.8][0]{$u_2$}
\psfrag{u0}[cc][][0.8][0]{$\iu$}
\psfrag{zeta}[cc][][0.8][0]{}
\psfrag{ud}[cc][][0.8][0]{$u_{\rm 1 d}$}
\includegraphics[scale=.8]{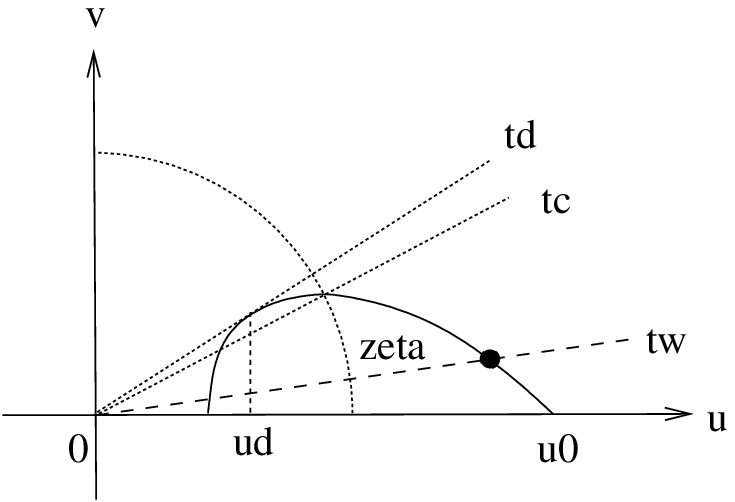}
\caption{The shock polar in the $\uu$--plane}\label{fig:polar}
\end{psfrags}
\end{figure}

Given $M_0>1$, $\rho_1$ and $\uu_1$ are determined by using the shock polar
in Fig. \ref{fig:polar} for steady potential flow.
For any wedge angle $\theta_{\rm w}\in (0,\theta_{\rm w}^{\rm s})$,
line $u_2=u_1\tan\theta_{\rm w}$ and the shock polar
intersect at a point $\uu_1$ with $|\uu_1| >c_1$ and $u_{11}<\iu$;
while, for any wedge angle $\theta_{\rm w}\in [\theta_{\rm w}^{\rm s}, \theta_{\rm w}^{\rm d})$,
they intersect at a point $\uu_1$ with $u_{11}>u_{\rm 1d}$
and $|\uu_1|<c_1$.
The intersection state $\uu_1$ is the velocity for steady potential flow
behind an oblique shock $\mathcal{S}_0$ attached to the wedge vertex with angle $\theta_{\rm w}$.
The strength of shock $\mathcal{S}_0$ is relatively weak compared to the other shock
given by the other intersection point on the shock polar, which is a
\emph{weak oblique shock},
and the corresponding state $\uu_1$ is a \emph{weak state}.

We also note that
states
$\uu_1$ depend smoothly on $\iu$ and $\theta_{\rm w}$, and
such states are supersonic when $\theta_{\rm w}\in (0,\theta_{\rm w}^{\rm s})$ and
subsonic when $\theta_{\rm w}\in [\theta_{\rm w}^s, \theta_{\rm w}^{\rm d})$.

Once $\uu_1$ is determined,
by \eqref{1-m} and \eqref{1-i},
the pseudo-potentials $\vphi_1$ and $\vphi_2$ below the weak oblique shock $\mathcal{S}_0$
and the normal shock $\mathcal{S}_1$
are respectively in the form of
\begin{equation}
\label{1-n}
\vphi_1=-\frac 12|\xxi|^2+ \uu_1\cdot\xxi,\qquad
\vphi_2=-\frac 12|\xxi|^2+ \uu_2\cdot \xxi+ k_2
\end{equation}
for constant states $\uu_1$ and $\uu_2$, and constant $k_2$. Then it follows
from \eqref{1-o} and \eqref{1-n} that the corresponding densities $\rho_1$ and $\rho_2$
below $S_0$ and $S_1$ are constants, respectively.
In particular, we have
\begin{equation}
\label{2-n1}
\rho_1^{\gam-1}=\rho_0^{\gamma-1}+\frac{\gam-1}{2}\big(\iu^2-|\uu_1|^2\big).
\end{equation}

Then {\bf Problem \ref{problem-5.2}} can be reformulated into the following
free boundary problem.

\smallskip
\begin{problem}[Free Boundary Problem]\label{fbp-a}
For $\theta_{\rm w}\in (0, \theta_{\rm w}^{\rm d})$,
find a free boundary {\rm (}curved shock{\rm )} $\shock$ and a function $\vphi$ defined in domain $\Om$,
as shown in Figs. {\rm \ref{fig:global-structure-1}}--{\rm \ref{fig:global-structure-2}},
such that $\vphi$ satisfies
\begin{itemize}
\item[\rm (i)]
Equation \eqref{2-1} in $\Om${\rm ;}
\item[\rm (ii)]
$\vphi=\ivphi$ and $\rho D\vphi\cdot\mathbf{n}_{\rm s}=\rho_0 D\ivphi\cdot\mathbf{n}_{\rm s}$ {on} $\shock${\rm ;}
\item[\rm (iii)]
$\vphi=\hat{\vphi}$ and $D\vphi=D\hat{\vphi}$ {on} $\Gamma_{\rm sonic}^1\cup\Gamma_{\rm sonic}^2\,\,$
when $\theta_{\rm w}\in (0, \theta_{\rm w}^s)$
and on $\Gamma_{\rm sonic}^2\cup \{O\}\,\,$ when $\theta_{\rm w}\in [\theta_{\rm w}^{\rm s}, \theta_{\rm w}^{\rm d})$
for $\hat{\vphi}:=\max(\vphi_1, \vphi_2)${\rm ;}
\item[\rm (iv)]
$D\vphi\cdot \mathbf{n}=0$ {on} $\Wedge$,
\end{itemize}
where $\mathbf{n}_{\rm s}$ and $\mathbf{n}$ are the interior unit normals to $\Omega$ on $\shock$ and $\Wedge$, respectively.
\end{problem}

Let $\vphi$ be a solution of {\bf Problem \ref{fbp-a}}
such that $\shock$ is a $C^1$--curve up to its endpoints and $\vphi\in C^1(\overline\Om)$.
To obtain a solution of {\bf Problem \ref{problem-5.2}} from $\vphi$,
we have two cases:

For $\theta_{\rm w}\in (0, \theta_{\rm w}^{\rm s})$,
we divide the half-plane $\{\xi_2\ge 0\}$ into five separate regions.
Let $\Om_{\mathcal{S}}$ be the unbounded domain below
curve $\ol{\mathcal{S}_0\cup\shock\cup \mathcal{S}_1}$
and above $\Wedge$
(see Fig. \ref{fig:global-structure-1}).
In $\Om_{\mathcal{S}}$, let $\Om_1$ be the bounded open domain enclosed
by $\mathcal{S}_0, \Gamma^1_{\rm sonic}$,
and $\{\xi_2=0\}$.
Set $\Om_{2}:=\Om_S\setminus \ol{(\Om_1\cup\Om)}$.
Define a function $\vphi_*$ in $\{\xi_2\ge 0\}$ by
\begin{equation}\label{extsol}
\vphi_*=
\begin{cases}
\ivphi& \qquad \text{in}\,  \Lambda\cap\{\xi_2\ge 0\}\setminus \Om_{\mathcal{S}},\\
\vphi_1& \qquad \text{in}\,\Om_1,\\
\vphi& \qquad \text{in}\, \Gamma^1_{\rm sonic}\cup\Om\cup\Gamma^2_{\rm sonic},\\
\vphi_2&\qquad \text{in}\,\Om_2.
\end{cases}
\end{equation}
By \eqref{1-i} and (iii) of {\bf Problem \ref{fbp-a}}, $\vphi_*$ is continuous in $\{\xi_2\ge 0\}\setminus\Om_S$
and $C^1$ in $\overline{\Om_{\mathcal{S}}}$.
In particular, $\vphi_*$ is $C^1$ across $\Gamma^1_{\rm sonic}\cup\Gamma^2_{\rm sonic}$.
Moreover, using (i)--(iii) of {\bf Problem \ref{fbp-a}}, we obtain
that  $\vphi_*$ is a global entropy solution of equation \eqref{2-1}
in $\Lambda\cap \{\xi_2>0\}$.

For $\theta_w\in [\theta_{\rm w}^{\rm s}, \theta_{\rm w}^{\rm d})$,
region $\Omega_1\cup \Gamma^1_{\rm sonic}$ in $\varphi_*$ reduces to one point $\{O\}$,
and the corresponding $\varphi_*$ is a global entropy solution of equation \eqref{2-1}
in $\Lambda\cap \{\xi_2>0\}$.

\smallskip
The first rigorous unsteady analysis of the steady supersonic weak shock solution
as the long-time behavior of an unsteady flow
is due to Elling-Liu \cite{EL2},
in which they
succeeded in establishing a stability theorem for an important class of physical
parameters determined by certain assumptions for the wedge angle $\theta_{\rm w}$
less than the sonic angle $\theta^{\rm s}_{\rm w}\in (0, \theta^{\rm d}_{\rm w})$
for potential flow.

Recently,  in Bae-Chen-Feldman \cite{BCF-14,BCF-16},
we have successfully remove the assumptions in Elling-Liu's theorem \cite{EL2}
and established the stability theorem for the steady (supersonic or transonic)
weak shock solutions as
the long-time asymptotics of
the global Prandtl-Meyer configurations for unsteady potential flow
for all the admissible physical parameters
even up to the detachment
angle $\theta^{\rm d}_{\rm w}$ (beyond the sonic angle $\theta^{\rm s}_{\rm w}<\theta^{\rm d}_{\rm w}$).
The global Prandtl-Meyer configurations involve two types of transonic transition --
discontinuous and continuous hyperbolic-elliptic phase transitions
for the fluid fields (transonic shocks and sonic circles).
To establish this theorem, we have first solved
the free boundary problem ({\bf Problem  \ref{fbp-a}}),
involving transonic shocks,
for all wedge angles
$\theta_{\rm w}\in (0, \theta_{\rm w}^{\rm d})$
by employing the
new techniques developed in Chen-Feldman
\cite{CF-Book} to obtain the monotonicity properties and uniform
{\it a priori} estimates for admissible solutions.
Therefore, we have achieved the existence of a self-similar weak solution
with higher regularity to {\bf Problem 5.1} for
all wedge angles
$\theta_{\rm w}$ up to the detachment angle $\theta_{\rm w}^{\rm d}$.

\smallskip
More precisely, to solve this free boundary problem,
we have followed the approach introduced in
Chen-Feldman \cite{CF-Book}.
We have first defined a class of admissible solutions  $\varphi$,
which are the solutions of Prandtl-Meyer configuration,
such that,
when $\theta_{\rm w}\in (0, \theta_{\rm w}^{\rm s})$,
equation (\ref{2-1}) is strictly elliptic for $\varphi$
in $\overline\Omega\setminus (\Gamma_{\rm sonic}^{1}\cup\Gamma_{\rm sonic}^{2})$,
$\max\{\varphi_{1}, \varphi_{2}\}\le\varphi\le \varphi_0$ holds in $\Omega$,
and the following monotonicity properties hold:
\begin{equation}\label{MonotoneProperty}
D(\varphi_0-\varphi)\cdot \mathbf{e}_{\mathcal{S}_1}\ge 0,
\quad D(\varphi_0-\varphi)\cdot \mathbf{e}_{\mathcal{S}_0}\le 0 \qquad\quad \mbox{in $\Omega$},
\end{equation}
where $\mathbf{e}_{\mathcal{S}_0}$ and $\mathbf{e}_{\mathcal{S}_{1}}$ are the unit tangential
directions to lines $\mathcal{S}_0$ and $\mathcal{S}_1$, respectively,
pointing to the positive $\xi_1$-direction.
The monotonicity properties in \eqref{MonotoneProperty} are the key to
ensure that the shock is a Lipschitz graph in a cone of directions,
so that the geometry of the problem is fixed,
among other consequences.
For $\theta_{\rm w}\in [\theta_{\rm w}^{\rm s}, \theta_{\rm w}^{\rm d})$,
admissible solutions have been defined similarly, with
corresponding changes to the structure of subsonic reflection solutions.

Another key step for solving {\bf Problem 5.3}
is to derive uniform {\it a priori} estimates for admissible solutions
for any wedge angle $\theta_{\rm w} \in [0, \theta_{\rm w}^{\rm d}-\varepsilon]$
for each $\varepsilon>0$.
In particular,
for fixed $\gam > 1$, $\iu>0$, and $\varepsilon>0$,
it has been proved that there exists a constant $C>0$ depending only on $\gam, \iu$, and $\varepsilon>0$
such that, for any $\theta\in(0,\theta^{\rm d}_{\rm w}-\varepsilon]$, a corresponding admissible
solution $\vphi$
satisfies
\begin{equation*}
{\rm dist}(\shock, B_{c_0}(\iu, 0))\ge C^{-1}>0.
\end{equation*}
This inequality plays an essential role to achieve the ellipticity of equation
\eqref{2-1}
in $\Om$. Once the ellipticity is achieved,
then we can obtain various apriori estimates of $\vphi$,
so that the Leray-Schauder degree argument can be employed to obtain
the existence for each $\theta_{\rm w} \in [0, \theta_{\rm w}^{\rm d}-\varepsilon]$
in the class of admissible solutions,
starting from the unique normal solution for $\theta_{\rm w}=0$.
Since $\varepsilon>0$ is arbitrary,
the existence of a weak solution for any $\theta_{\rm w}\in(0,\theta_{\rm w}^{\rm d})$
can be established.

More details can be found in Bae-Chen-Feldman \cite{BCF-14,BCF-16}; see also
Chen-Feldman \cite{CF-Book}.

\smallskip
More recently, in Chen-Feldman-Xiang \cite{ChenFeldmanXiang},
we have also established the strict convexity of the curved transonic
part of the free boundary in
the Prandtl-Meyer configurations described above.
In order to prove the convexity, we employ the global properties of
admissible solutions, including the existence of the cone of monotonicity
discussed above.

The existence results in Bae-Chen-Feldman \cite{BCF-14,BCF-16}
indicate that the steady weak supersonic/transonic shock solutions are the asymptotic limits
of the dynamic self-similar solutions, the Prandtl-Meyer configurations,
in the sense of \eqref{1-k-c} in {\bf Problem 5.1}.

On the other hand, it is shown in Elling \cite{Elling2} and  Bae-Chen-Feldman \cite{BCF-16} that,
for each $\gam >1$, there is no self-similar {\it strong} Prandtl-Meyer configuration
for the unsteady potential flow in the class of admissible solutions ({\it cf.} \cite{BCF-16}).
This means that the situation for the dynamic stability of the strong steady oblique shocks
is more sensitive.

\section{Further Problems and Remarks}

In \S 1--\S 5, we have surveyed some recent developments on the static stability
of the weak and strong steady shock solutions for the wedge problem;
we have also presented the
recent results on the dynamic stability of the weak steady supersonic/transonic
shock solutions for potential flow.
These indicate that the weak supersonic/transonic oblique shocks are both stable,
and it is more sensitive for the dynamic stability of the steady strong transonic shocks,
which require further mathematical understanding.
Moreover, there are many other open problems in this direction,
which require further investigations.

\smallskip
When the deviation of vorticity become significant, the full Euler equations are required.
It is still open how the Prandtl-Meyer configurations can be constructed
for the full Euler flow.
As seen in \S 5,
we have understood  the mathematical
difficulties relatively well for the transonic shocks,
the Kelydsh degeneracy near the sonic arcs,
and the corner
between the transonic shock and the sonic arcs
for the nonlinear second-order elliptic equations,
as well as a one-point singularity at the wedge vertex between the attached shock and
the wedge boundary for the transition of state $(1)$ from the supersonic to subsonic
states
when the wedge angle increases across the sonic angle
up to the detachment angle.
On the other hand,
when the flow is pseudo-subsonic, the system consists of two
transport-type equations and two nonlinear equations of
mixed hyperbolic-elliptic type.
Therefore, in general, the full Euler system
is of {\em composite-mixed hyperbolic-elliptic type};
see Chen-Feldman \cite{CF-Book}.
Then the following two new features
for this problem for the isentropic and/or full Euler
equations still need to be understood:

\smallskip
(i) Solutions of transport-type equations with
rough coefficients and stationary transport velocity;

(ii)  Estimates of the vorticity of the pseudo-velocity.

\smallskip
\noindent
Indeed, a similar calculation as in Serre \cite{Serre1} has shown difficulties
in estimating the vorticity.
It is possible that the vorticity has some singularities in the
region, perhaps near the wedge boundary and/or corner.
In fact, even for potential flow, for the wedge
angle $\theta_{\rm w}\in (\theta_{\rm w}^{\rm s}, \theta_{\rm w}^{\rm d})$,
the second derivatives of the
velocity potential, {\it i.e.}, the first derivatives of the velocity,
may blow up at the wedge corner.

\smallskip
For the global stability of three or higher dimensional (M-D) transonic shocks
in steady supersonic flow past M-D wedges, the situation
is much more sensitive than that for the $2$-D case, which requires
more careful rigorous mathematical analysis.
In Chen-Fang \cite{CFang2},  we developed a nonlinear approach and employed it to
establish the stability of weak shock solutions containing
a transonic shock for potential flow with respect to the
M-D perturbation of the wedge boundary in appropriate
function spaces.
To achieve this,
we first formulated the stability problem as an M-D free boundary problem
for nonlinear elliptic equations.
Then we introduced the partial
hodograph transformation to reduce the free boundary problem into
a fixed boundary value problem near a background solution
with fully nonlinear boundary conditions
for second-order nonlinear elliptic equations
in an unbounded domain in M-D.
To solve this reduced problem,
we linearized the nonlinear problem on the background shock solution and
then, after solving this linearized elliptic problem,
we developed a nonlinear iteration scheme that was proved to be contractive,
which implies the convergence of the scheme to yield the desired results.
It would be interesting to investigate further problems
for
the stability of M-D shocks
in steady supersonic flow past M-D wedges.
In this regard, we notice that an instability result has been observed
in Liang-Xu-Yin \cite{LiangXuYin}.

\smallskip
Conical flow ({\it i.e.}, cylindrically symmetric flow
with respect to an axis) occurs in many physical
situations. For instance, it occurs at the conical nose of a
projectile facing a supersonic stream of air ({\it cf.} \cite{CourantF}).
The global stability of conical supersonic
shocks has been studied
in Liu-Lien
\cite{LL}
in the class of $BV$ solutions when the cone vertex angle is small,
and Chen \cite{Sxchen2} and Chen-Xin-Yin \cite{CXY} in the class of
smooth solutions away from the conical shock when the
perturbed cone is sufficiently close to the straight-sided cone.
The stability of transonic shocks in 3-D steady
flow past a perturbed cone had been a longstanding open problem.
For the 2-D wedge case, the equations do not involve
such singular terms, and the flow past the straight-sided wedge is
piecewise constant. However, for the 3-D conical case,
the governing equations have a singularity at the cone vertex and
the flow past the straight-sided cone is self-similar, but no
longer piecewise constant.
These cause additional difficulties for the stability problem.
In Chen-Fang \cite{CFang1}, we developed techniques to handle
the singular terms in
the equations and the singularity of the solutions.
Our main results indicate that the self-similar transonic
shock is conditionally stable with respect to the conical
perturbation of the cone boundary and the upstream flow in
appropriate function spaces.
That is, it was proved that the transonic shock and downstream
flow in our solutions are close to the unperturbed self-similar
transonic shock and downstream flow under the conical
perturbation, and
the slope
of the shock asymptotically tends to the slope of the
unperturbed self-similar shock at infinity.
These results were obtained by first formulating the stability
problem as a free boundary problem and introducing
a coordinate
transformation to reduce the free boundary problem into a fixed
boundary value problem for a singular nonlinear elliptic system.
Then we
developed an iteration scheme that consists of two iteration mappings:
One is for an iteration of approximate transonic shocks, and
the other is for an iteration of the corresponding boundary value
problems for the singular nonlinear systems for given approximate
shocks. To ensure the well-definedness and contraction
property of the iteration mappings, it is essential to establish the
well-posedness for a corresponding singular linearized elliptic
equation, especially the stability with respect to the coefficients
of the equation,  and to obtain the estimates of its solutions
reflecting their singularity at the cone vertex and decay at
infinity. The approach is to employ key features of the equation,
to introduce appropriate solution spaces, and to apply a Fredholm-type
theorem in Maz'ya-Plamenevski\v{\i} \cite{MP} to establish the
existence of solutions by showing the uniqueness in the solution
spaces.

\smallskip
Another important direction is to analyze the detached shocks off the wedge
when the
supersonic
flow onto the wedge whose angle is larger than the detachment angle.

\smallskip
Finally, many fundamental problems in this direction
are still wide open, and their solution requires further new techniques,
approaches and ideas,
which deserve our special attention.

\bigskip

\Acknowledgements{
The materials presented in this paper contain direct and indirect
contributions of my collaborators
Myoungjean Bae, Jun Chen, Beixiang Fang, Mikhail Feldman, Tianhong Li,
Wei Xiang, Dianwen Zhu, and Yongqian Zhang.
The work of Gui-Qiang G. Chen was supported in part by
the US National Science Foundation under Grants
DMS-0935967 and DMS-0807551,
the UK
Engineering and Physical Sciences Research Council under Grants
EP/E035027/1 and EP/L015811/1,
the National Natural Science Foundation of China (under joint project Grant 10728101),
and the Royal Society--Wolfson Research Merit Award (UK).}


\end{document}